\documentclass[letter,12pt]{article}
\usepackage{epsfig,psfrag,amsmath,amssymb,latexsym,eufrak,amscd,amsfonts,graphicx,fancyhdr}
\usepackage[margin=2cm,includeheadfoot]{geometry}
\pdfoutput=1

\newtheorem{theorem}{Theorem}[section]

\newtheorem{corollary}{Corollary}[section]

\newtheorem{example}{Example}[section]
\newtheorem{lemma}{Lemma}[section]
\newtheorem{proposition}{Proposition}[section]

\newcommand{\eop}{\mbox{ \vrule height7pt width7pt depth0pt}}

\numberwithin{equation}{section}

\title{Random modification effect in the size of the fluctuation of the LCS of two sequences of i.i.d. blocks}
\author{Heinrich Matzinger\footnote{School of Mathematics, Georgia Institute of Technology, 686 Cherry Street, GA 30332-0160 Atlanta, USA}\quad Felipe Torres\footnote{Fakult\"at f\"ur Mathematik, Universit\"at Bielefeld, Postfach 100131 D-33501 Bielefeld, Germany}}
\date{\today}

\begin{document}
\maketitle
\thispagestyle{empty}
\abstract{\noindent The problem of the order of the fluctuation of the Longest Common Subsequence (LCS) of two independent sequences has been open for decades. There exist contradicting conjectures on the topic, \cite{CS} and \cite{W4}. In the present article, we consider a special model of i.i.d. sequences made out of blocks. A {\it block} is a contiguous substring consisting only of one type of symbol. Our model allows only three possible block lengths, each been equiprobable picked up. For i.i.d. sequences with equiprobable symbols, the blocks are independent of each other. For this model, we introduce a random operation (random modification) on the blocks of one of the sequences. In this article, for our block model, we show the techniques to prove the following: if we suppose that the random modification increases the length of the LCS with high probability, then the order of the fluctuation of the LCS is as conjectured by Waterman \cite{W4}. This result is a key technical part in the study of the size of the fluctuation of the LCS for sequences of i.i.d. blocks, developed in \cite{MT1}.}


\section{Model and main results}
In general trough this paper, $X$ and $Y$ will denoted two finite strings over a finite alphabet $\Sigma$. A common subsequence of $X$ and $Y$ is a subsequence which is a subsequence of $X$ as well as of $Y$. A Longest Common Subsequence of $X$ and $Y$ (denoted simply by {\rm LCS} of $X$ and $Y$, or only LCS when the context is clear enough) is a common subsequence of $X$ and $Y$ of maximal length. For a motivation on why to study the LCS problem, the reader can look at \cite{MT1,FelipePhD}.
\\
\\
Let $l>0$ be an integer parameter. Let $B_{X1},B_{X2},\ldots$ and $B_{Y1},B_{Y2},\ldots$ be two i.i.d. sequences independent of each other such that:
\begin{eqnarray}
\mathrm{P}(B_{Xi}=l-1)=\mathrm{P}(B_{Xi}=l)=\mathrm{P}(B_{Xi}=l+1)&=&1/3 \nonumber\\
\mathrm{P}(B_{Yi}=l-1)=\mathrm{P}(B_{Yi}=l)=\mathrm{P}(B_{Yi}=l+1)&=&1/3. \nonumber
\end{eqnarray}
We call the runs of $0$'s and $1$'s blocks. Let $X^{\infty}=X_1X_2 X_3\ldots$ be the binary sequence so that the $i$-th block has length $B_{Xi}$ where $X_1$ is choosen $0$ with probability $1/2$ or $1$ with probability $1/2$. Similarly let $Y^{\infty}=Y_1Y_2 Y_3\ldots$ be the binary sequence so that the $i$-th block has length $B_{Yi}$ and $Y_1$ is choosen $0$ with probability $1/2$ or $1$ with probability $1/2$.

\begin{example}
Assume that $X_1=0$ and  $B_{X1}=3$, $B_{X2}=4$ and $B_{X3}=2$. Then we have  that the sequence  $X^\infty$ starts as follows $X^\infty=000111100\cdots$ meaning that in $X^\infty$ the first block consists of three 0's, the second block consists of four 1's, the third block consists of two 0's, etc.
\end{example}
Let $X$ denote the sequence obtained by only taking the first $n$ bits of $X^\infty$, namely $X=X_1X_2X_3\ldots X_n$ and similarly $Y=Y_1Y_2Y_3\ldots Y_n.$ Let $L_n$ denote the length of the {\rm LCS} of $X$ and $Y$, $L_n:=|{\rm LCS}(X,Y)|.$ 
\\
The main result of \cite{MT1,FelipePhD} states that for $l$ large enough, the order of the fluctuation  of $L_n$ is $n$:
\begin{theorem} \label{maintheorem}
There exists $l_0$ so that for all $l\geq l_0$ we have that:
\[\mathrm{VAR}[L_n]=\Theta(n)\]
for $n$ large enough.
\end{theorem}
In \cite{MT1,FelipePhD} the authors showed that theorem \ref{maintheorem} is equivalent to proving that ``a certain random modification has a biased effect on $L_n$''. This is a technique with similar approches in other papers (for instance see \cite{JJMM}, \cite{BM}). So the main difficulty is actually proving that the random modification has typically a biased effect on the {\rm LCS}, which for the block model is connected to a constrained optimization problem \cite{MT1,FelipePhD}. This random modification is performed as follows: we choose at random in $X$ a block of length $l-1$ and at random one block of length $l+1$, this means that all the blocks in $X$ of length $l-1$ have the same probability to be chosen and then we pick one of those blocks of length $l-1$ up and also that all the blocks in $X$ of length $l+1$ have the same probability to be chosen and we pick one of those blocks of length $l+1$ up. Then we change the length of both these blocks to $l$. The resulting new sequence is denoted by $\tilde{X}$. Let $\tilde{L}_n$ denote the length of the {\rm LCS} after our modification of $X$. Hence: 
\[\tilde{L}_n:=|{\rm LCS}(\tilde{X},Y)|.\]
If we can prove that our block length changing operation has typically a biased effect on the {\rm LCS} than the order of the fluctuation of $L_n$ is $\sqrt {n}$. This is the content of the next theorem:
\begin{theorem} \label{theorem2}
Assume that there exists $\epsilon>0$ and $\alpha>0$ not depending on $n$ such that for all $n$ large enough we have:
\begin{equation}\label{fundamental}
\mathrm{P}\left(\;\;  \mathrm{E}[\tilde{L}_n-L_n|X,Y]\geq \epsilon\;\;\right)\geq 1-\exp(-n^\alpha).
\end{equation}
Then,
$$\mathrm{VAR}[L_n]=\Theta(n)$$
for $n$ large enough.
\end{theorem}
The above theorem reduces the problem of the order of fluctuation to proving that our random modification has typically a higher probability to lead to an increase than to a decrease in score. {\bf The main result of this article is theorem \ref{theorem2}}.
\\
\\
\noindent A very useful tool we often use is the Azuma-Hoeffding theorem. The following is a version of it for martingales (for a proof see \cite{GS}):
\begin{theorem}(Hoeffding's inequality)\label{hoeffding}
Let $(V,\mathfrak{F})$ be a martingale, and suppose that there exists a sequence $\mathfrak{a}_1,\mathfrak{a}_2,\cdots$ of real numbers such that
\[ \mathrm{P}(|V_n-V_{n-1}| \le \mathfrak{a}_n)=1\]
for all $n$. Then:
\begin{equation}
\mathrm{P}(|V_n-V_0| \ge v) \le 2 \exp\Big\{ -\frac{1}{2} v^2 \Big/ \sum_{i=1}^n \mathfrak{a}^2_i \Big\}\label{hoeffdingIneq}
\end{equation}
for every $v>0$.
\end{theorem}
We also will use a corollary of the above theorem, for some intermediate bounds:
\begin{corollary}\label{hoeffdingco}
Let $a>0$ be constant and $V_1,V_2,\dots$ be an i.i.d sequence of random bounded variables such that:
\[ {\rm P}(|V_i-{\rm E}[V_i]| \le a)=1\]
for every $i=1,2,\dots$ Then for every $\Delta>0$, we have that:
\begin{equation}\label{corAH}
{\rm P}\left(\,\left|\frac{V_1+\dots+V_n}{n}-{\rm E}[V_1] \right|\ge \Delta\right) \le 2\exp\left( -\frac{\Delta^2}{2a^2}\cdot n\right)
\end{equation}
\end{corollary}

\section{Random modification effect in the fluctuation} \label{secRM}
We are going to prove theorem \ref{theorem2} which states that ${\rm VAR}[L_n]=\Theta(n)$ holds if there exist $\epsilon, \alpha>0$ not depending on $n$ such that:
\begin{equation}\label{tildeepsilon}
\mathrm{P}\left(\;\;  \mathrm{E}[\tilde{L}_n-L_n|X,Y]\geq \epsilon\;\;\right)\geq
1-\exp(-n^\alpha).
\end{equation} 
for all $n$ large enough. We have omitted some of the proofs for shortness reasons, but all the details can be looked at \cite{FelipePhD}.
\\
\\
Note that if $\mathcal{Z}$ is a random variable with ${\rm VAR}[\mathcal{Z}]=\Theta(n)$ and $f$ is a map which tends to increase linearly, then for $\mathcal{W}=f(\mathcal{Z})$, we also have the order ${\rm VAR}[\mathcal{W}]=\Theta(n)$. The map $f$ can be even a random map but must be independent of $\mathcal{Z}$. The exact basic result (\cite{BM}, lemma 3.2) goes as follows:

\begin{lemma}\label{bonetto}
Let $c>0$ be a constant. Assume that $g:\mathbb{R}\rightarrow\mathbb{R}$ is a map which is everywhere differentiable and such that for all $x\in\mathbb{R}$ we have:
\[\frac{dg(x)}{dx}\ge c.\]
Let $B$ be a random variable such that $\mathrm{E}[|g(B)|] < +\infty$. Then:
$${\rm VAR}[g(B)]\geq c^2 \cdot {\rm VAR}[B].$$
\end{lemma}
In the present context, we need a slightly different version:
\begin{lemma} \label{bonettorefined}
Let $\epsilon,m>0$ be constants and $f: \mathbb{Z} \to \mathbb{Z}$ be a map such that for all $z_1 \le z_2$ the following two conditions hold:
\begin{equation}\label{superlog}
z_2 - z_1 \ge m \Rightarrow f(z_2)-f(z_1) \ge \frac{\epsilon}{8}(z_2-z_1)
\end{equation}
\begin{equation}\label{sublog}
\exists\,\,\beta>0:\; z_2 - z_1 < m \Rightarrow f(z_2)-f(z_1) \le \beta (z_2-z_1)
\end{equation}
Let $B$ be a random variable such that $\mathrm{E}[|f(B)|]\le+\infty$. Then:
\begin{equation}
{\rm VAR}[(f(B)] \ge \frac{\epsilon^2}{64}\left( 1-16\frac{(\epsilon/8+\beta)m}{\epsilon\sqrt{{\rm VAR}[B]}} \right){\rm VAR}[B]
\end{equation}
\end{lemma}
\vspace{12pt}{\bf Proof.}
Let $h:\mathbb{Z} \to \mathbb{Z}$ be a map defined from $f$ as follows: for a given $z \in \mathbb{Z}$ choose $k \ge 2$ such that $z \in [km,(k+1)m]$ and compute
\[ h(z) = \left(\frac{f((k+1)m)-f(km)}{m}\right)(z-km)+f(km)\]
then $h(z)$ is just the linear interpolation of $f(z)$ in $[km,(k+1)m]$. It is easy to see that $h$ satisfies the conditions of lemma \ref{bonetto} for $c=\epsilon/8$. Then:
\begin{equation}\label{inter1}
{\rm VAR}[h(B)] \ge \frac{\epsilon^2}{64}{\rm VAR}[B]
\end{equation}
We want to estimate the distance between the random variables $h(B)$ and $f(B)$. First, we note that from \ref{superlog} and by the definition of $h$, the following inequalities hold for $km \le B \le (k+1)m$:
\[\frac{\epsilon}{8}(B-km)+f(km) \le f(B),h(B) \le \frac{\epsilon}{8}(B-(k+1)m)+f((k+1)m)\]
looking at conditions \ref{superlog}, \ref{sublog} and the inequalities above we get
\begin{eqnarray}
|h(B)-f(B)| &\le& | \frac{\epsilon}{8}(B-km)+f(km) - \frac{\epsilon}{8}(B-(k+1)m)+f((k+1)m)| \nonumber \\
&\le& \frac{\epsilon}{8}m+|f((k+1)m)-f(km)|  \nonumber \\
&\le& \left( \frac{\epsilon}{8}+\beta\right)m \nonumber
\end{eqnarray}
and by using the last inequality above:
\begin{equation}\label{inter2}
{\rm VAR}[f(B)-h(B)] \le \left( \frac{\epsilon}{8}+\beta\right)^2 m^2.
\end{equation}
Since $f(B)=h(B)+(f(B)-h(B))$ we can apply triangular inequality and find:
\[ \sqrt{{\rm VAR}[f(B)]} \ge \sqrt{{\rm VAR}[h(B)]}-\sqrt{{\rm VAR}[f(B)-h(B)]},\]
hence we have:
\begin{eqnarray}
{\rm VAR}[(f(B)] &\ge& {\rm VAR}[h(B)] - 2\sqrt{{\rm VAR}[h(B)]} \cdot \sqrt{{\rm VAR}[f(B)-h(B)]} \nonumber \\
&=& {\rm VAR}[h(B)]\left( 1-2\frac{\sqrt{{\rm VAR}[f(B)-h(B)]}}{\sqrt{{\rm VAR}[h(B)]}}\right) \nonumber
\end{eqnarray}
Finally, applying the inequalities \ref{inter1} and \ref{inter2} to the last inequality above, we get:
\[{\rm VAR}[(f(B)] \ge \frac{\epsilon^2}{64}\left( 1-16\frac{(\epsilon/8+\beta)m}{\epsilon\sqrt{{\rm VAR}[B]}} \right){\rm VAR}[B].\quad\eop\]

\vspace{12pt}
\noindent Hence to prove that ${\rm VAR}[L_n]=\Theta(n)$, we try to represent $L_n$
as $f(\mathcal{Z})$ where $f$ is a random map which tends to increase
linearly on a certain scale and $\mathcal{Z}$ is a random variable having
fluctuation of order $\sqrt{n}$.
\subsection{Random modifications and the variables $(T,Z,R)$}\label{subTZ}
Let $N_l$ denote the number of blocks in $X$ of length $l$,
whilst $N_{l-1}$, resp. $N_{l+1}$ denote the number of blocks of 
length $l-1$, resp $l+1$
in $X$. Let us define the following three random variables:
\begin{eqnarray}
T&:=&N_l+N_{l-1}+N_{l+1}\label{T}\\
Z&:=&N_l-N_{l-1}-N_{l+1}\label{Z} \\
R&:=&n-(\,l\, N_l+(l+1)\, N_{l+1}+(l-1)\, N_{l-1}\,)\label{R}
\end{eqnarray}
Note that when we know the values of $(T,Z,R)$ we can determine the values of $N_{l-1},N_l$ and $N_{l+1}$ as a linear function by using the definitions of $T,Z$ and $R$ as follows:
\begin{equation}\label{linearTZ}
\left(\begin{array}{c} N_{l-1}(T,Z,R) \\ N_l(T,Z,R) \\ N_{l+1}(T,Z,R) \end{array}\right)=
\left(\begin{array}{cc} (2l+1)/4 & -1/4 \\ 1/2 & 1/2 \\ -(2l-1)/4 & -1/4\end{array}\right) \left( \begin{array}{c} T \\ Z \end{array}\right)+ \left(\begin{array}{c} -(n-R)/2 \\ 0 \\ (n-R)/2 \end{array} \right)
\end{equation}
The variable $R$ represents what is left in $X$ after the last block of length $l-1, l$ or $l+1$.
\begin{example}\label{ejemploR}
Let us consider the sequence $X=000111100011001$ for $l=3$ and $n=15$. We see that $N_{l-1}=2$, $N_l=2$ and $N_{l+1}=1$, hence $T=5,Z=-1$ and $R=1$. Also, the block $1$ at the end of $X$ has length strictly smaller than $l-1$ which also means that $R=1$. In this case is easy to interpret what $R$ is since the last block in $X$ has length strictly less than $l-1$. Let us see a different situation. Let us take again $l=3$ and now consider $B_{X1}=2,B_{X2}=3,B_{X3}=4,B_{X4}=3,B_{X5}=2,B_{X6}=4,\dots$ such that $X^\infty=001110000111001111\cdots$ Take $n=16$ so that $X=0011100001110011$. Here the last block of $X$ has length $l-1=2$ which should imply (using the point of view of the last situation) that $R=0$. But, notice that the block in $X^\infty$ corresponding to $B_{X6}$ was cut when we took $X$. In this case, we say that the last block in $X$ corresponds to the rest so $R=2$ and therefore $N_{l-1}=2,N_l=2$ and $N_{l+1}=1$, then $T=5$ and $Z=-1$. We take this convention on $R$, even if the definition \ref{R} is not the exact one, because of the simplifications later during the computation of the joint distribution of $N_{l-1},N_l,N_{l+1}$.
\end{example}
Let us roughly explain the main idea behind this subsection. Assume that we have a random couple $(V,W)$ which can take on a finite number of values only. We also assume the joint distribution $\mathcal{L}(V,W)$ to be given. To simulate $(V,W)$, we could first simulate $V$ using the marginal law $\mathcal{L}(V)$. We would obtain a numeric value $v_0$. Then, we could simulate $W$ using the conditional law $\mathcal{L}(W|V=v_0)$ and obtain the numeric value $w_0$. The couple $(v_0,w_0)$ has joint distribution $\mathcal{L}(V,W)$. Another less efficient possibility is to simulate for each (non-random) value $v$ that $V$ can take, a value for $W$ with distribution
$\mathcal{L}(W|V=v)$. Call the numeric value $w(v)$. Then, we would simulate $V$ with distribution $\mathcal{L}(V)$ and obtain a numeric value $v_0$. Then, for $W$ we would take among all the values which we have si\-mu\-lated, the one corresponding to $V=v_0$. In this manner, we get $(v_0,w(v_0))$. This couple has the distribution $\mathcal{L}(V,W)$ and this does not even require that
we simulate the  different $w(v)$'s independently of each other. Only, $V$ needs to be simulated independently of the assignment $v\mapsto w(v)$.
\\
\\
\noindent We are going to do the above simulation scheme with $V$ being $(T,Z,R)$ and $W$ being the rest of the information in $(X,Y)$. More precisely, for all possible $(t,z,r)$ non-random values, we simulate $X$ conditional on $(T,Z,R)=(t,z,r)$. The resulting string is denoted by $X_{(t,z,r)}$ and has thus distribution
$$\mathcal{L}(X_{(t,z,r)})=\mathcal{L}(X\;|\;(T,Z,R)=(t,z,r)\;).$$
Let $L_n(t,z,r)$ denote the length of the LCS 
$$L_n(t,z,r):=|{\rm LCS}(X_{(t,z,r)},Y)|.$$
We assume that the simulation of the string $X_{(t,z,r)}$ is done independently of $(T,R,Z)$ and of $Y$. In this manner, we get that $L_n(T,Z,R)$ has same distribution as $L_n=|{\rm LCS}(X,Y)|.$ So to prove that ${\rm VAR}[L_n]=\Theta(n)$, it
is enough to prove that
\begin{equation} \label{LnTR}
{\rm VAR}[L_n(T,Z,R)]=\Theta(n).
\end{equation}
We saw at the beginning of this section (see lemma \ref{bonetto} and \ref{bonettorefined}), that when we transform a variable having variance of order $\Theta(n)$ with a map which tends to increase linearly, then the resulting variable has variance of order $\Theta(n)$. It is easy to see that ${\rm VAR}[Z]=\Theta(n)$ (see also lemma \ref{VARZD}). Hence to prove \ref{LnTR}, it is enough to show that
with high probability the (random) map
$$z\mapsto L_n(T,z,R)$$ tends to increase linearly (on the appropriate scale and on a domain on which $Z$ typically takes its value). That means, we need to show that we can si\-mu\-late the values $L_n(t,z,r)$ in such a manner to get the desired distribution $\mathcal{L}(X|(T,Z,R)=(t,z,r))$ as well as the desired linear increase of the  map $z\mapsto L_n(T,z,R)$. This is achieved by simulating $X_{(t,z,r)}$ in the following way: for a given value $(t,r)$, so that ${\rm P}((T,R)=(t,r))\neq 0$, we take a left most (left most to be defined later)  value $z_0$ and simulate a string with distribution equal to the conditional distribution of $X$ given $(T,Z,R)=(t,z_0,r)$. That resulting string is denoted by $X_{(t,z_0,r)}$. Then, we apply the random modification to $X_{(t,z_0,r)}$. This means, we choose one block of length $l-1$ and one block of length $l+1$ at random in $X_{(t,z_0,r)}$ and turn them both into length $l$. The resulting string is denoted by $X_{(t,z_0+4,r)}$. Then, we choose at random in $X_{(t,z_0+4,r)}$ a block of lenght $l-1$ and a block of lenght $l+1$ and turn them both into length $l$. The new string which we obtain in this manner is denoted by  $X_{(t,z_0+8,r)}$. We keep repeating this same operation to obtain the sequence of strings 
\begin{equation} \label{sequence}
X_{(t,z_0,r)},X_{(t,z_0+4,r)},X_{(t,z_0+8,r)},\ldots.
\end{equation}
For each value of $(t,r)$ with ${\rm P}((T,R)=(t,r))\neq 0$ we obtain two finite sequences of strings: first \ref{sequence} and then
\[X_{(t,z_0+2,r)},X_{(t,z_0+6,r)},X_{(t,z_0+10,r)},\ldots.\]
by a similar procedure. Namely, after $X_{(t,z_0+2,r)}$ is generated with distribution
$X$ conditional on $(T,Z,R)=(t,z+2,R)$, the 
subsequent strings are obtained by applying sucessively
the random modification tilde, which chooses at random in the string a block of length $l-1$ and a block of length $l+1$ and turn them both into length $l$.
\\
\\
\noindent Recall that in this section we assume that our random modification has a biased effect of $\epsilon>0$ on the LCS, so that with high probability
$${\rm E}[\tilde{L}_n-L_n\;|\;X,Y\;]\geq \epsilon.$$
Hence, it follows that the map $z\mapsto L_n(T,z,R)$ tends with high probability to increase with slope close to $\epsilon$ on a constant time scale $\ln n$ (the constant must be taken large enough though, see lemma \ref{highA} and proposition \ref{driftepsilon}). In other words, since the random modification has a biased positive effect, the map $z\mapsto L_n(T,z,R)$ behaves like a random walk with drift $\epsilon$. The only thing which remains to be proved is that with our scheme of using the random modification, the strings $X_{(t,z,r)}$
have the right distribution, i.e. the distribution of $X$ conditional on $(T,Z,R)=(t,z,r)$. This is proved in lemma \ref{possZ}.
\\
\\
\noindent We have so far summarized the idea which explains why the biased effect of the random modification implies ${\rm VAR}[L_n]=\Theta(n)$. There is one more detail which we should mention and which makes notations a little more difficult. To prove that $z\mapsto L_n(T,z,R)$ tends to increase linearly we use the biased effect on the LCS for the random modification. However, this bias holds with high probability for $X$ and not for $X_{(t,z,r)}$. When we look at the conditional distribution of $X$ given $(T,Z,R)=(t,z,r)$, we divide by the probability
\begin{equation}\label{TZR}
{\rm P}((T,Z,R)=(t,z,r)).
\end{equation}
The string $X_{(t,z,r)}$ has distribution of $X$ conditional on $(T,Z,R)=(t,z,r)$. So for the biased effect to have large probability also for $X_{(t,z,r)}$ (and not just for $X$), we need the probability \ref{TZR} to not be too small. To assure this, we will restrict ourselves to ``typical'' values for $(T,Z,R)$. We will consider only values for $(T,Z)$ which lie in an interval $D=D_T\times D_Z$ (see definition below \ref{dominioD}) and prove that any possible value $(t,z)\in D_z\times D_t$ has polynomially bounded probability (see lemma \ref{localcentral}).
\\
\\
\noindent  Let us now give all the details:

\begin{proposition}\label{propCLTnl}
Given $\epsilon >0$ there exist constants $1\le k_1,k_2,k_3 \le k^*$ all not depending on $n$ but on $\epsilon$ such that:
\begin{equation}\label{CLTnl}
\mathrm{P}\left(\left| \frac{N_{l-1} -\frac{n}{3l}}{\sqrt{n}} \right| \le  k_1\right), \quad \mathrm{P}\left(\left| \frac{N_l -\frac{n}{3l}}{\sqrt{n}} \right| \le  k_2\right), \quad \mathrm{P}\left(\left| \frac{N_{l+1} -\frac{n}{3l}}{\sqrt{n}} \right|\le k_3\right) \ge 1-\epsilon
\end{equation}
for every $n$ large enough.
\end{proposition}

\vspace{12pt}
\noindent We will need later the following lemma:
\begin{lemma}
There exists $c>0$ not depending on $n$ such that:
\begin{equation} \label{0.9}
\mathrm{P}\left(T\in\left[\frac{n}{l}-c\sqrt{n},\frac{n}{l}+c\sqrt{n}\right], Z\in\left[-\frac{n}{3l}-c\sqrt{n},-\frac{n}{3l}+c\sqrt{n}\right]\right)\geq 0.9
\end{equation}
\end{lemma}



\vspace{12pt}
\noindent Let $D$ denote the domain
\begin{equation}\label{dominioD}
D:=\left[\frac{n}{l}-c\sqrt{n},\frac{n}{l}+c\sqrt{n}\right] \times \left[-\frac{n}{3l}-c\sqrt{n},-\frac{n}{3l}+c\sqrt{n}\right]
\end{equation}
and let 
\begin{eqnarray}
D_T&:=&\left[\frac{n}{l}-c\sqrt{n},\frac{n}{l}+c\sqrt{n}\right] \nonumber \\
D_Z&:=&\left[-\frac{n}{3l}-c\sqrt{n},-\frac{n}{3l}+c\sqrt{n}\right] \nonumber
\end{eqnarray}
hence,
$$D=D_T\times D_Z.$$
Given $(t,z) \in D$ such that $(T,Z)=(t,z)$ we have
\[N_{l-1}(t,z)+N_l(t,z)+N_{l+1}(t,z)=t.\]
The probability for a realization of $N_{l-1},N_l$ and $N_{l+1}$ is given by:
\begin{eqnarray}
{\rm P} (T=t,Z=z,R=r)&=&{N_{l-1}(t,z)+N_l(t,z)+N_{l+1}(t,z) \choose N_{l-1}(t,z)\;\;N_l(t,z)\;\;N_{l+1}(t,z)} \left(\frac{1}{3}\right)^t \cdot {\rm P}(B_{X1} > r) \nonumber \\
&=&\frac{t!}{(N_{l-1}(t,z))!\;(N_l(t,z))!\;(N_{l+1}(t,z))!} \left(\frac{1}{3}\right)^t \cdot {\rm P}(B_{X1} > r)\nonumber \\ \label{mult}
\end{eqnarray}
where the probability ${\rm P}(B_{X1}>r)={\rm P}(R=r)$ is due to the convention of $R$ described in the example \ref{ejemploR}. Finally, due to \ref{linearTZ}, for any $n_1,n_2,n_3 \in \mathbb{N}$ the conditional joint distribution
\[\mathrm{P}(N_{l-1}(T,Z,R)=n_1,N_l(T,Z,R)=n_2,N_{l+1}(T,Z,R)=n_3\,|\,R=r)\]
is multinomial.
\begin{lemma}\label{localcentral}
There exists $k_0>0$ not depending on $n$ (but depending on $c$) such that for every $(t,z)\in D$ and $r<l+1$ for which the probability ${\rm P}((T,Z,R)=(t,z,r)) \ne 0$, we have that:
$$\mathrm{P}((T,Z,R)=(t,z,r))\geq \frac{k_0}{n}$$
for every $n$ large enough.
\end{lemma}

\vspace{12pt}
\noindent Note that for any variables $X$ and $Y$ we have (see for example \cite{ROSS})
\begin{equation} \label{VARXY}
{\rm VAR}[Y]=\mathrm{E}[{\rm VAR}[Y|X]]+{\rm VAR}[\mathrm{E}[Y|X]]\geq \mathrm{E}[{\rm VAR}[Y|X]].
\end{equation}
Let $O$ be the random variable which is equal to one
when $(T,Z)$ is in $D$ and $0$ otherwise.\\
We can now use inequalities \ref{0.9} and \ref{VARXY} to find
\begin{equation}
\label{VAR09}
 {\rm VAR}[L_n]\geq \mathrm{E}[{\rm VAR}[L_n|O]]\geq {\rm VAR}[L_n|O=1]\cdot\mathrm{P}(O=1) \ge 0.9{\rm VAR}[L_n|O=1]
\end{equation}
Next for every $(t,z)$ in $D$ and $r<l+1$ we are going to simulate
the random variable $L_n$  conditional on
$(T,Z,R)=(t,z,r)$. We denote the result by $L_n(t,z,r)$.
In other words, the distribution
of $L_n(t,z)$ is equal to
$$\mathcal{L}(L_n(t,z,r))=
\mathcal{L}(L_n|(T,Z,R)=(t,z,r)).$$
Let $(T_D,Z_D)$ denote a variable having the distribution
of $(T,Z)$ conditional on the event 
$(T,Z)\in D$.
We assume that all the $L_n(t,z,r)$ are independent of $(T_D,Z_D)$.
Then, we get that
$$L_n(T_D,Z_D,R)$$
has same distribution as $L_n$ conditional on $(T,Z)\in D$.
Hence, we get
\begin{equation}
\label{VAR1}
{\rm VAR}[L_n|O=1]={\rm VAR}[L_n(T_D,Z_D,R)]
\end{equation}
By using \ref{VARXY}, we find
\begin{equation} \label{VAR2}
{\rm VAR}[L_n(T_D,Z_D,R)]]\geq
\mathrm{E}[\;{\rm VAR}[L_n(T_D,Z_D,R)|T_D,R\,]\;].
\end{equation}
Note that for $L_n(T_D,Z_D,R)$ to have the same distribution as $L_n$ conditional on $(T,Z)\in D$ and on $R=r$, the variables $L_n(t,z,r)$ do not need to be independent of each other. We are next going to explain how we simulate the variables $L_n(t,z,r)$ a bit more in detail as before. We simulate a string $X_{(t,z,r)}$ having the distribution of the string $X$ conditional on the event $(T,Z,R)=(t,z,r)$. Then we put
$$L_n(t,z,r)=|{\rm LCS}(X_{(t,z,r)},Y)|.$$
Next, let us describe how we simulate $X_{(t,z,r)}$ based on what was roughly explained at the beginning of subsection \ref{subTZ}. Given $t_0\in D_T$ the most left element in $D_T$ and $r_0<l-1$, we are going to simulate $X_{(t_0,z,r_0)}$ for $z\in D_Z$ only if ${\rm P}((T,Z,R)=(t_0,z,r_0)) \ne 0$.  We simulate $X_{(t_0,z_0,r_0)}$ so that it has distribution $\mathcal{L}(\,X|(T,Z,R)=(t_0,z_0,r_0)\,)$. Next, we simulate $X_{(t_0,z_0+2,r_0)}$ by choosing in $X$, with the same probability, a block of length $l-1$ either a block of length $l+1$ and change its length to $l$. The next realization we simulate is $X_{(t_0,z_0+4,r_0)}$ by choosing in $X$, with the same probability, a block of length $l-1$ and a block of length $l+1$ and change their lengths to $l$ (this is our usual random modification). Then by induction we simulate 
$$\{X_{(t_0,z_0+4i,r_0)}:i=1,2,\dots\}$$
with our usual random modification and later
$$\{X_{(t_0,z_0+2+4i,r_0)}:i=1,2,\dots\}$$
just starting with $X_{(t_0,z_0+2,r_0)}$ and performing our usual random modification to get each $X_{(t_0,z_0+6,r_0)},X_{(t_0,z_0+10,r_0)}, X_{(t_0,z_0+14,r_0)}$, etc. Both inductions run untill indexes $i_0$, resp. $i_0^*$, satisfying:
\begin{eqnarray}
z_0+4i_0 &\le& -\frac{n}{3l}+c\sqrt{n} \quad\Rightarrow\quad i_0 \le \sqrt{n}\nonumber \\
z_0+2+4i_0^* &\le&-\frac{n}{3l}+c\sqrt{n} \quad\Rightarrow\quad i_0^* \le \frac{\sqrt{n}-1}{2}\nonumber
\end{eqnarray}
For simplicity, let us call $z_0,z_1=z_0+2,z_2=z_0+4,\dots,z_d$ all the values which $Z$ takes. After we have simulated $X_{(t_0,z_0,r_0)},X_{(t_0,z_1,r_0)},\dots,X_{(t_0,z_d,r_0)}$ we fix $t_1=t_0+1$ and repeat all the procedure again starting with the simulation of $X_{(t_1,z_0,r_0)}$. We keep taking $t_2 < t_3 < t_4\dots$ all natural numbers in $D_T$ to finish all the simulation of $\{X_{(t,z,r_0)}:t\in D_T, \,z=z_0,z_1,\dots,z_d\}$. Once we have finished with that, we take $r_1<l-1$ natural number and do all the simulation above starting with $X_{(t_0,z_0,r_1)}$ only if ${\rm P}((t_0,z_0,r_1)) \ne 0$. Finally, we obtain the complete sequence $\{X_{(t,z,r)}:t\in D_T, \,z=z_0,z_1,\dots,z_d\,,\,r=0,\dots,l-2\}$, where each $(t,z,r)$ has probability ${\rm P}((T,Z,R)=(t,z,r)) \ne 0$.

\vspace{12pt}
\noindent We need to verify that this operation give us the equiprobable distribution. This is the content of the next lemma:
\begin{lemma}\label{possZ}
Assume that $X_{(t,z,r)}$ is distributed according to 
$$\mathcal{L}(X|(T,Z,R)=(t,z,r)).$$
Choose at random (with equal probability) in the string $X_{(t,z,r)}$ a block of length $l+1$ and $l-1$ and modify them to have both length $l$. Then the resulting string has distribution
$$\mathcal{L}(X|(T,Z,R)=(t,z+4,r)).$$
\end{lemma}
\vspace{12pt}{\bf Proof.}
Because of our linear equation system \ref{linearTZ}, we have that conditioning on $T,Z,R$ is equivalent to conditioning on $(N_{l-1},N_l,N_{l+1})$. As mentioned, $X_{(t,z,r)}$ denotes a string of length $n$, having the distribution of $X$ conditional on
$(T,Z,R)=(t,z,r)$. We denote by $\tilde{X}_{(t,z,r)}$ the string we obtain by performing our random modification on $X_{(t,z,r)}$. In other words, $\tilde{X}_{(t,z,r)}$ is obtained by choosing a block of length $l+1$ and a block
of length $l-1$ at random in $X_{(t,z,r)}$ and changing them both to length $l$. Let $(n_1,n_2,n_3)$ be the number of
blocks of length $l-1$, $l$ and $l+1$ corresponding to $(t,z,r)$. In other words, $n_1$, $n_2$
and $n_3$ are given by the linear system of equation \ref{linearTZ} when $N_{l-1}=n_1,N_l=n_2,N_{l+1}=n_3$ and $T=t,Z=z,R=r$. We have
$${\rm P}(N_1=n_1,N_2=n_2,N_3=n_3|T=t,Z=z,R=r)=1.$$
The distribution of the random string  $X_{(t,z,r)}$ is the uniform distribution on $\xi^n(t,z,r)$. Here, $\xi^n(t,z,r)$ denotes the set of  strings of length $n$, which consists only of blocks of length $l-1$, $l$ and $l+1$, such that the total number of blocks is $t$, whilst the number of blocks of length $l$ minus the number of blocks of length $l-1$ and $l+1$ is $z$. We also request that the rest block at the end has length $r$.  We can describe $\xi^n(t,z,r)$ equivalently as the set of all strings consisting exactly of $n_1$ blocks of length $l-1$, $n_2$ blocks of length $l$ and $n_3$ blocks of length $l+1$, no other blocks allowed except a rest block at the end which has length strictly less than $l-1$.
In other words, the random string $X_{(t,z,r)}$ is such that the number of blocks of length $l-1$, $l$ and $l+1$ is determined, only the order in which these blocks appear varies. Among others, each possible realization for $X_{(t,z,r)}$ which has non-zero probability has the same probability:
\begin{equation}
\label{binom}
\binom{n_1+n_2+n_3}{n_1\;n_2\;n_3}^{-1}
\end{equation}
When we apply the random modification, the variable $T$ stays the same,
the variable $Z$ increases by $4$ and the variable $R$ stays the same.
\\
\\
\noindent Since the distribution of $X$ conditional on $(T,Z,R)$ is the uniform distribution on the appropriate set of strings, we have the following: for proving that $\tilde{X}_{(t,z,r)}$ has distribution of $X$ conditional on $(T,Z,R)=(t,z+4,r)$ it is enough to show that
its distribution is  the uniform distribution on $\xi^n(t,z+4,r)$. For this, let $\tilde{x}$ denote a (non-random) element of $\xi^n(t,z+4,r)$. Hence, the number of blocks in $\tilde{x}$ of length $l-1$, $l$, resp $l+1$ is $n_1-1$, $n_2+2$, resp. $n_3-1$. The probability 
$${\rm P}(\tilde{X}_{(t,z,r)}=\tilde{x})$$
can be calculated as follows: if we only know $\tilde{x}$, any  block of length $l$ of $\tilde{x}$ could be the block which had lenght $l-1$ and has been turned into length $l$ by the {\it tilde operation} (choosing blocks at random and changing their lenghts). Same thing for the block which had length $l+1$. But when we know these two blocks, then the string before the random modification is uniquely determined. Let $x$ be such a string which could lead to $\tilde{x}$ after the random modification. There are hence $\tilde{n}_2\cdot (\tilde{n}_2-1)$ such strings (here, $\tilde{n}_2=n_2+2$, so that $\tilde{n}_2$ denotes the number of blocks of length $l$ in $\tilde{x}$).
The probability, given  $X_{(t,z,r)}=x$, that the random string turns out to be $\tilde{x}$ is equal to $1/(n_1\cdot n_3)$. As a matter of fact,
among the $n_1$ blocks of lenght $l-1$, there is exactly one which needs to be randomly modified. Similarly, among the $n_3$ blocks of lenght $l+1$, there is exactly one which needs to be changed into length $l$ in order to obtain the string $\tilde{x}$. Hence,
\begin{equation} \label{n1n3}
{\rm P}(\tilde{X}_{(t,r,z)}=\tilde{x}| X_{(t,z,r)}=x)=\frac{1}{n_1\cdot n_3}.
\end{equation}
Let $\xi^{n*}$ denote the set of all strings which could lead to $\tilde{x}$ if we apply the random modification to them. We saw that there are  $(n_2+2)(n_2+1)$ elements in the set $\xi^{n*}$. By law of total probability, we have
{\small \begin{equation} \label{Idonotknow}
{\rm P}(\tilde{X}_{(t,z,r)}=\tilde{x})= \sum_{x\in \xi^{n*}} {\rm P}(\tilde{X}=\tilde{x}|X=x){\rm P}(X_{(t,z,r)}=x)
=\sum_{x\in \xi^{n*}} \frac{1}{n_1\cdot n_3} \binom{n_1+n_2+n_3}{n_1\;n_2\;n_3}^{-1}
\end{equation}}
The last equation above was obtained using \ref{n1n3} and \ref{binom}. Note that  the sum on the most right of equation in \ref{Idonotknow}, is a sum of $(n_2+2)(n_2+1)$ equal terms. This leads to
$${\rm P}(\tilde{X}_{(t,z,r)}=\tilde{x})= \frac{(n_2+2)(n_2+1)}{n_1\cdot n_3}\binom{n_1+n_2+n_3}{n_1\;n_2\;n_3}^{-1}.$$
The formula on the right side above  does not depend on $\tilde{x}$. Hence, this proves that $\tilde{X}_{(t,z,r)}$ has the uniform distribution on the set of strings $\xi^n(t,z+4,r)$. But the uniform distribution is the distribution of $X$ conditional on $(T,Z,R)=(t,z+4,r)$. That is, we have proven that
$$\mathcal{L}(\tilde{X}_{(t,r,z)})=\mathcal{L}(X|(T,Z,R)=(t,z+4,r)),$$
which finishes this proof. $\quad\eop$

\vspace{12pt}
\noindent Note that we have seen what happens with the variables $T,Z,R$ after our random modification, let us see what happens with the length of the {\rm LCS} after our random modification. In what follows, we always consider a triplet of values $(t,z,r)$ such that ${\rm P}((T,Z,R)=(t,z,r)) \ne 0$. For any $\epsilon>0$ let $U^n_{t,r}(\epsilon)$ denote the event that the map
$$D_Z\rightarrow \mathbb{N} \; : \; z\mapsto L_n(t,z,r)$$
is increasing with a slope of at least $\epsilon/8$ on a scale
$c_2\ln(n)$ where $c_2>0$ is a large constant not depending on $n$. 
More precisely, $U_{t,r}^n(\epsilon)$ is the event that for any 
$z_1,z_2$ in $D_Z$,
with $z_2-z_1\geq c_2 \ln(n)$ we have
$$L_n(t,z_2,r)-L_n(t,z_1,r)\geq (z_2-z_1)\epsilon/8.$$
The event $U^n_{t,r}(\epsilon)$ has large probability because
we assumed that inequality \ref{tildeepsilon} holds.
Hence $z\mapsto L_n(t,z,r)$ can be viewed somehow as behaving
like a random walk with drift $\epsilon$. In the next lemma we will show this looking at the event $U^n(\epsilon)$:
\[U^n(\epsilon):= \bigcap_{t\in D_T,\,\, r <l+1}U_{t,r}^n(\epsilon).\]
\begin{lemma} \label{highA}
Given $\epsilon>0$, take $\alpha$ from inequality \ref{tildeepsilon} (theorem \ref{theorem2}) and $c_2$ to be big enough but not depending on $n$, for example $c_2 \ge \frac{80}{\epsilon^2}$ depending on $\epsilon$. Then, there exists a constant $k_*>0$ not depending on $n$ but on $\alpha$ and on $c_2$ such that:
\begin{equation}\label{expsmalldrift}
\mathrm{P}(U^{nc}(\epsilon))\leq \frac{k_*}{n^2}
\end{equation}
for $n$ large enough, provided \ref{tildeepsilon} holds. 
\end{lemma}
\vspace{12pt}{\bf Proof.}
We are going to define an event $\mathcal{U}(\epsilon)$ for any $\epsilon>0$.  Let $\mathcal{U}_{(t,z,r)}(\epsilon)$ be the event that the expected conditional increase is larger than $\epsilon$ when we introduce the random change into $X_{(t,z,r)}$. More precisely, let
$\mathcal{U}_{(t,z,r)}^n(\epsilon)$ be the event that
\begin{equation}
\label{epsilon4Ltz}
\mathrm{E}[\;L_n(t,z+4,r)-L_n(t,z,r)|X_{(t,z,r)},Y\;]\geq \epsilon
\end{equation}
Let
$$\mathcal{U}^n(\epsilon):=\bigcap_{(t,z)\in D,\,\,r<l+1}\mathcal{U}_{(t,z,r)}^n(\epsilon).$$ 
hence
\begin{equation}
\label{calUnc}\mathrm{P}(\mathcal{U}^{nc}(\epsilon)) \le \sum_{(t,z)\in D,\,\,r<l+1}\mathrm{P}(\mathcal{U}_{(t,z,r)}^{nc}(\epsilon)).
\end{equation}
Note that inequality \ref{tildeepsilon} provides a bound for the probability that
the conditional expected increase of {\rm LCS} due to our random
modification not being larger or equal to $\epsilon$. That
probability bound is $\exp(-n^\alpha)$. The only problem
is that the bound is for $X$ and $Y$ whilst
the event $\mathcal{U}^n_{(t,z,r)}(\epsilon)$ is for $X_{(t,z,r)}$ and $Y$.
By going on to conditional probability we must multiply the 
probability by $\mathrm{P}((T,Z,R)=(t,z,r))$.
Hence we find
\begin{equation}
\label{in3}
\mathrm{P}(\mathcal{U}_{(t,z,r)}^{nc}(\epsilon))\leq \frac{\exp(-n^\alpha)}{\mathrm{P}((T,Z,R)=(t,z,r))}.
\end{equation}
We can next use the lower bound on $\mathrm{P}((T,Z,R)=(t,z,r))$
provided by lemma \ref{localcentral} for all values $(t,z)\in D$ and $r<l+1$
to inequality \ref{in3} and obtain
\begin{equation}
\label{in4}
\mathrm{P}(\mathcal{U}_{(t,z,r)}^{nc}(\epsilon))\leq \frac{1}{k_0}\cdot n\cdot\exp(-n^\alpha).
\end{equation}
which still gives an exponentially small bound in $n$.
Applying now \ref{in4} to inequality \ref{calUnc}, we obtain
\begin{equation}
\label{Unc2}\mathrm{P}(\mathcal{U}^{nc}(\epsilon))\leq \frac{4lc^2}{k_0}\cdot n^2\cdot\exp(-n^\alpha).
\end{equation}
Which is an exponentially small bound in $n$. Note that when the event $\mathcal{U}^n(\epsilon)$ holds, we have
that $z\mapsto L_n(t,z,r)$ behaves like a random walk
with drift $\epsilon$. Let us formalize this. As before, let $\{z_0,z_1,z_2,\dots,z_d\}$ be the set for the admissible values of $Z$.
For fixed $t\in D_T$ and $r<l+1$, we are going to define
$L_n^*(t,z)$ inductively for $z\in \{z_0,z_1,z_2,\dots,z_d\}$. Let us define $L_n^*(t,z,r):=L_n(t,z,r)$ for every $z\in\{z_0,z_1,z_2,\dots,z_d\}$. Given $\tilde{z}\in \{z_0,z_1,z_2,\dots,z_d-4\}$ let us define $L_n^*(t,\tilde{z}+4,r)$ as follows:
\[L_n^*(t,\tilde{z}+4,r) = \left\{\begin{array}{cl}
L_n(t,\tilde{z}+4,r) & \mbox{ if\,\, $\mathcal{U}^n_{(t,s,r)}(\epsilon)$ hold for all $s\in\{z_0,z_1,\dots,\tilde{z}\}$} \\
L^*_n(t,\tilde{z},r)+\epsilon &\mbox{ otherwise}
\end{array}\right.\]
Note that when the event $\mathcal{U}^n(\epsilon)$ holds,
then $L_n(t,z,r)$ and $L_n^*(t,z,r)$
are identical for all $t\in D_T$, $r<l+1$ and $z\in\{z_0,z_1,\dots,z_d\}$.
Let $\mathcal{V}_{t,r}^{n}(\epsilon)$ be the event that the map
$$D_Z\rightarrow \mathbb{N}\; : \; z\mapsto L_n^*(t,z,r)$$
is increasing with a slope of at least $\epsilon/8$ on a scale
$c_2\ln n$.

\vspace{12pt}
\noindent Let $\mathcal{V}^{n}(\epsilon)$ be the event $$\mathcal{V}^{n}(\epsilon):= \bigcap_{t\in D_T,\,\,r<l+1}\mathcal{V}_{t,r}^{n}(\epsilon).$$
Hence by using proposition \ref{driftepsilon} we have that:
\begin{equation} \label{RWdrift}
\mathrm{P}(\mathcal{V}^{nc}(\epsilon))\leq \sum_{t\in D_T,\,\,r<l+1} \mathrm{P}(\mathcal{V}_t^{nc}(\epsilon))\le \sum_{t\in D_T,\,\,r<l+1} 2n^{-\tau} \le 4lc\,n^{0.5-\tau}
\end{equation}
where $\tau=\frac{\epsilon^2\,c_2}{32}$. When $\mathcal{U}^n(\epsilon)$ holds then $\mathcal{V}^{n}(\epsilon)$ and $U^{n}(\epsilon)$ are equivalent. Hence
$$\mathcal{U}^n(\epsilon)\cap \mathcal{V}^{n}(\epsilon)\subset U^n(\epsilon)$$
Hence by using \ref{Unc2} and \ref{RWdrift} we get:
\begin{equation}\label{bb0}
\mathrm{P}(U^{nc}(\epsilon))\leq \mathrm{P}(\mathcal{U}^{nc}(\epsilon))+\mathrm{P}(\mathcal{V}^{nc}(\epsilon)) \le \frac{4lc^2}{k_0}\cdot n^2\cdot\exp(-n^\alpha) + 4lc\,n^{0.5-\tau}
\end{equation}
To show that the last inequality gives us a rate of convergence to zero as a constant divided by a polynomial in $n$, we try now to get a closed form for the inequality supposing extra information for the involved constants. 

\vspace{12pt}
\noindent Taking $c_2 \ge \frac{80}{\epsilon^2}$ we have the following bound for the exponent:
\[ 0.5-\tau \le-2\]
therefore we can bound
\begin{equation}\label{bb1}
4lc\,n^{0.5-\tau} \le \frac{4lc}{n^2}\,.
\end{equation}
Also, we have that:
\begin{equation}
n^2 \exp\left(n^{-\alpha}\right) \le \frac{1}{n^2}\label{bb2}
\end{equation}
holds for $n$ large enough. So, by using \ref{bb1} and \ref{bb2} in \ref{bb0} we can finally bound:
\begin{eqnarray}
\mathrm{P}(U^{nc}(\epsilon))\leq \mathrm{P}(\mathcal{U}^{nc}(\epsilon))+\mathrm{P}(\mathcal{V}^{nc}(\epsilon)) &\le& 4lc^2\tilde{c}_2n^2\cdot\exp(-n^\alpha) + 4lc\,n^{0.5-\tau} \nonumber \\
&\le& (4lc^2\tilde{c}_2+4lc) \cdot \frac{1}{n^2}\nonumber
\end{eqnarray}
for $n$ large enough, which ends the proof with $k_*=4lc^2k_0+4lc$.
$\quad\eop$

\begin{proposition}\label{driftepsilon}
Given $\epsilon>0$, let $\mathcal{V}^{n}_{t,r}(\epsilon)$ denote the event that the map $z \mapsto L^*(t,z,r)$ is increasing with a slope at least $\epsilon/8$ on a scale $c_2 \ln(n)$. Given $t \in D_T$, $r<l+1$ and $z_1,z_2 \in D_Z$ such that $z_2-z_1 \ge c_2\ln(n)$ we have the following inequality:
\[{\rm P}\left( \mathcal{V}^{nc}_{t,r}(\epsilon) \right) \le 2\,n^{-\tau}\]
where $\,\tau=\frac{\epsilon^2\,c_2}{32}$.
\end{proposition}
\vspace{12pt}{\bf Proof.}
Let $z_1,z_2\in D_Z$ such that $z_1 < z_2$. In order to simplify the notation, let us assume that $z_2-z_1$ can be dived by 4 and denote $\frac{z_2-z_1}{4}=m\in \mathbb{N}$. Let $z_0$ be the most left point of $D_Z$. Given $\epsilon>0$, let us remember that $\mathcal{V}^n_{t,r}(\epsilon)$ is the event such that the following inequality holds:
\[ L^*(t,z_2,r)-L^*(t,z_1,r) \ge \frac{\epsilon}{8}.\]
Now let us define the filtration $\mathfrak{F}_0 \subset \mathfrak{F}_1 \subset \cdots \subset \mathfrak{F}_m$ as follows:
\[ \mathfrak{F}_i:=\sigma\left( X_{(t,z_0,r)},X_{(t,z_1,r)},\dots,X_{(t,z_1+4i,r)}\,;\,Y\right)\]
for $i=1,\dots,m$.  Let us denote
\[ e_i=E[\,L^*_n(t,z_1+4(i+1),r)- L^*_n(t,z_1+4i,r)\,|\,\mathfrak{F}_i\,]\]
and define a martingale $M_0,M_1,\dots,M_m$ with respect to the filtration $\mathfrak{F}_0 \subset \mathfrak{F}_1 \subset \cdots \subset \mathfrak{F}_m$ as follows:
\begin{eqnarray}
M_0 &:=&L^*_n(t,z_1,r) \nonumber \\
M_{i+1} -M_i&:=&L^*_n(t,z_1+4(i+1),r)- L^*_n(t,z_1+4i,r)-e_i \nonumber
\end{eqnarray}
for $i=1,\dots,m$. By definition of the map $z \mapsto L^*_n(t,z,r)$ we have an expected increase of at least $\epsilon$ every time $z$ gets increased by 4, so that the expected increase of 
\[E[\,L^*_n(t,z_1+4(i+1),r)- L^*_n(t,z_1+4i,r)\,]\]
is at least $\epsilon$ which implies that the following inequality
\begin{equation}
e_i \ge \epsilon \label{incepmart}
\end{equation}
is satisfied almost surely for every $0=1,\dots,m$. We can write the increase of the map $z \mapsto L^*_n(t,z)$ in terms of the martingale $M_0,\dots,M_m$ in the following way:
\begin{equation}\label{mart1}
 L^*(t,z_2,r)-L^*(t,z_1,r) = M_m-M_0+\sum_{i=0}^{m-1} e_i
\end{equation}
Now, we are ready to estimate the probability of $\mathcal{V}^{nc}_{t,r}(\epsilon)$:
\begin{eqnarray}
{\rm P}\left( \mathcal{V}^{nc}_{t,r}(\epsilon) \right) &=&{\rm P}\left(  L^*(t,z_2,r)-L^*(t,z_1,r) \le \frac{\epsilon}{8}(z_2-z_1)\right) \nonumber \\
\mbox{(by equality \ref{mart1})}&\le&{\rm P}\left( M_m-M_0+\sum_{i=0}^{m-1} e_i \le \frac{\epsilon}{8}(z_2-z_1)\right) \nonumber \\
&=&{\rm P}\left( M_m-M_0 \le \frac{\epsilon}{8}(z_2-z_1)-\sum_{i=0}^{m-1} e_i\right) \nonumber \\
\mbox{(by \ref{incepmart} and $\,z_2-z_1=4m$)}&\le&{\rm P}\left( M_m-M_0 \le \frac{\epsilon}{8}(z_2-z_1)-\frac{\epsilon}{4}(z_2-z_1)\right) \nonumber \\
&=&{\rm P}\left( M_m-M_0 \le -\frac{\epsilon}{8}(z_2-z_1)\right) \label{mart2}
\end{eqnarray}
At this point we want to use Azuma-Hoeffding inequality \ref{hoeffding}. For this, we note that for every $i=1,\dots,m$ we have
\[{\rm P}(|M_{i+1}-M_i| \le 1)=1\]
since $\epsilon<1$ and we take $v=\frac{\epsilon}{8}(z_2-z_1)$ for writing down:
\begin{eqnarray}
{\rm P}\left( M_m-M_0 \le -\frac{\epsilon}{8}(z_2-z_1)\right) &\le& 2\exp\left( -\frac{v^2}{2m}\right) \nonumber \\
\mbox{(by using $\,z_2-z_1=4m$)}&=& 2\exp\left( -\frac{\epsilon^2}{32}(z_2-z_1)\right) \label{mart3}
\end{eqnarray}
Combining together \ref{mart2} and \ref{mart3} we finally have:
\[ {\rm P}\left( \mathcal{V}^{nc}_{t,r}(\epsilon) \right) \le 2\exp\left( -\frac{\epsilon^2}{32}(z_2-z_1)\right)\]
from where, after taking $z_2-z_1 \le c_2\ln(n)$, we have:
\[{\rm P}\left( \mathcal{V}^{nc}_{t,r}(\epsilon) \right) \le 2\exp\left( -\frac{\epsilon^2\,c_2}{32}\ln(n)\right)=2\,n^{-\frac{\epsilon^2\,c_2}{32}}\]
which finishes the proof
$\quad\eop$

\vspace{24pt}
\noindent Note that by law of total probability $\mathrm{E}[\;{\rm VAR}[L_n(T_D,Z_D,R)|T_D,R\,]\;]$
is equal to :
{\small $$\mathrm{P}(U^n(\epsilon))\mathrm{E}[\;{\rm VAR}[L_n(T_D,Z_D,R)|T_D,R\,]\;|U^n(\epsilon)]
+\mathrm{P}(U^{nc}(\epsilon))\mathrm{E}[\;{\rm VAR}[L_n(T_D,Z_D,R)|T_D,R\,]\;|U^{nc}(\epsilon)],$$}
for every $\epsilon>0$ and hence:
\begin{equation}
\label{VARLnT}
\mathrm{E}[\;{\rm VAR}[L_n(T_D,Z_D,R)|T_D,R\,]\;]
\geq \mathrm{P}(U^n(\epsilon))\mathrm{E}[\;{\rm VAR}[L_n(T_D,Z_D,R)|T_D,R\,]\;|U^n(\epsilon)]
\end{equation}
Now, conditional on the event $U^n(\epsilon)$ holding,
we have that the random map:
$$D_Z\rightarrow\mathbb{N}\;:\; z\mapsto L_n(t,z,r)$$
has a slope of at least $\epsilon/8$ on a scale of $c_2\ln(n)$ (as in proposition \ref{driftepsilon}) for any $t\in D_T$ and $r<l+1$, then:
\begin{eqnarray}
z_2-z_1 \ge c_2\ln(n) &\Rightarrow& L_n(t,z_2,r)-L_n(t,z_1,r) \le \frac{\epsilon}{8}\,(z_2-z_1)\nonumber\\
z_2-z_1<c_2\ln(n) &\Rightarrow& L_n(t,z_2,r)-L_n(t,z_1,r) \le 2\,(z_2-z_1)\nonumber
\end{eqnarray}
hold. Hence, conditional on $U^n(\epsilon)$, we can apply  lemma \ref{bonettorefined} and obtain:
{\footnotesize \begin{equation}\label{VAR|A}
{\rm VAR}[L_n(t,Z_D,R)|T_D=t,R=r,U^n(\epsilon)]\geq \frac{\epsilon^2}{64}\left( 1-16\frac{(\epsilon/8+2)c_2 \ln(n)}{\epsilon\sqrt{{\rm VAR}[Z_D|T_D=t,R=r]}} \right){\rm VAR}[Z_D|T_D=t,R=r]
\end{equation}}
The next results give us an uniform bound for ${\rm VAR}[Z_D|T_D=t,R=r]$ for all $t\in D_T$.
\begin{lemma}
There exists a constant $K>0$ not depending on $n$ such that:
\begin{equation}\label{BinScope}
1-\frac{K}{\sqrt{n}} \le \frac{P (Z_D=z+4|T_D=t,R=r)}{P (Z_D=z|T_D=t,R=r)} \le 1+\frac{K}{\sqrt{n}}
\end{equation}
for every $(t,z) \in D$, $r<l+1$ and $n$ large enough.
\end{lemma}
\begin{lemma}\label{VARZD}
There exists a constant $C>0$ not depending on $n$ such that:
\[ {\rm VAR}[Z_D|T_D=t,R=r] \ge C \cdot n\]
for every $t\in D_T$, $r<l+1$ and for every $n$ large enough.
\end{lemma}

\vspace{12pt}
\noindent Using the bound in lemma \ref{VARZD} we get the following inequality:
\begin{equation}
\left( 1-16\frac{(\epsilon/8+2)c_2 \ln(n)}{\epsilon\sqrt{{\rm VAR}[Z_D|T_D=t,R=r]}} \right)\ge \left( 1-16\frac{(\epsilon/8+2)c_2 }{\epsilon\sqrt{C}}\cdot\frac{\ln(n)}{\sqrt{n}} \right)\geq 0.5\label{above1}
\end{equation}
for $n$ large enough. Using inequality \ref{above1} above with inequality \ref{VAR|A} we find:
\[{\rm VAR}[L_n(t,Z_D,R)|T_D=t,R=r,U^n(\epsilon)] \geq \frac{\epsilon^2}{64}\,\,0.5 \cdot{\rm VAR}[Z_D|T_D=t,R=r].\]
Using again lemma \ref{VARZD} we find that the left side of the above inequality is larger than $\frac{C\epsilon^2}{128}n$ and hence:
\begin{equation}
\label{last2}
\mathrm{E}[\;{\rm VAR}[L_n(T_D,Z_D,R)|T_D,R]\;|U^n(\epsilon)\;]\geq \frac{C\epsilon^2}{128}n
\end{equation}
We can now combine inequalities \ref{VAR09},
\ref{VAR1}, \ref{VAR2}, \ref{VARLnT} and \ref{last2}
to obtain:
$${\rm VAR}[L_n]\geq \mathrm{P}(U^n(\epsilon))\frac{C\epsilon^2}{1000}n$$
and plugging in the lower bound for $\mathrm{P}(U^n(\epsilon))$ obtained in \ref{expsmalldrift} (lemma \ref{highA}) we get:
$${\rm VAR}[L_n]\geq \frac{C\epsilon^2}{1000}\,n\left(1-\frac{k_*}{n^2}\right)$$
with $k_*>0$ is the constant from lemma \ref{highA}. This expression is a lower bound of order $\Theta(n)$ for ${\rm VAR}[L_n]$. Hence, we have finished proving the statement of the result in theorem \ref{theorem2}.

\vspace{36pt}
\noindent {\bf \large Acknowledgments}
\vspace{12pt}

\noindent The authors would like to thank the support of the German Science Foundation (DFG) through the International Graduate College "Stochastics and Real World Models" (IRTG 1132) at Bielefeld University and  through the Collavorative Research Center 701 "Spectral Structures and Topological Methods in Mathematics" (CRC 701) at Bielefeld University.


\begin{thebibliography}{9}
\bibitem{CS} V. Chvatal and D. Sankoff. {\it Longest common subsequences of two random sequences}. J. Appl. Probability, 12 : 306--315, 1975.
\bibitem{W4} M. S. Waterman. {\it Estimating statistical significance of sequence alignments}. Phil. Trans. R. Soc. Lond. B, 344:383-390, 1994.
\bibitem{MT1} H. Matzinger, Torres, F. {\it Fluctuation of the longest common subsequence for sequences of independent blocks}. Submitted, 2010.
\bibitem{FelipePhD} Torres, F. {\it On the probabilistic longest common subsequence problem for sequences of independent blocks}. Ph.D. thesis, University of Bielefeld. March 2009. Online http://bieson.ub.uni-bielefeld.de/volltexte/2009/1473/
\bibitem{JJMM} J. Lember, H. Matzinger. {\it Standard Deviation of the Longest Common Subsequence}. Ann. Probab. Volume 37, Number 3: 1192-1235, 2009.
\bibitem{BM} F. Bonetto and H. Matzinger. {\it Fluctuations of the longest common subsequence in the case of 2- and 3-letter alphabets}. Latin American Journal of Probability and Mathematics, Volume 2:195--216, 2006.
\bibitem{GS} G. Grimmett. and D. Strizaker. {\it Probability and Random Processes}, Oxford University Press, 2001. Third edition.
\bibitem{ROSS} S.M. Ross, \textit{Introduction Probability Models}. Academic Press, 8 edition, 2002.




\end{thebibliography}
\end{document}